\documentclass{article}
\usepackage{amssymb,amsmath}
\usepackage{graphics}

\def\qed{\hfill {\hbox{${\vcenter{\vbox{               
   \hrule height 0.4pt\hbox{\vrule width 0.4pt height 6pt
   \kern5pt\vrule width 0.4pt}\hrule height 0.4pt}}}$}}}
  
\def\tr{\triangleright}
\def\tl{\triangleleft}
  
\newtheorem{theorem}{Theorem}
\newtheorem{definition}{Definition}
\newtheorem{lemma}[theorem]{Lemma}
\newtheorem{proposition}[theorem]{Proposition}
\newtheorem{corollary}[theorem]{Corollary}
\newtheorem{example}{Example}

\newenvironment{proof}[1][Proof]{\smallskip\noindent{\bf #1.}\quad}%
{\qed\par\medskip}

\author{{\begin{tabular}{c} Gabriel Murillo \\
\small{\texttt{gmuri002@student.ucr.edu}}\end{tabular}}
\and
{\begin{tabular}{c} Sam Nelson \\
\small{\texttt{knots@esotericka.org}}\end{tabular}}
\and {\begin{tabular}{c} Anthony Thompson\\
\small{\texttt{athom005@student.ucr.edu}}\end{tabular}}
\and
\small{University of California, Riverside, 900 University Avenue,
Riverside, CA, 92521 }}

\date{}

\title{\Large \textbf{Matrices and finite Alexander quandles}}

\begin{document}

\maketitle

\begin{abstract}
We study the question of whether a finite quandle specified by a matrix is 
isomorphic to an Alexander quandle. We first collect some standard observations
about necessary conditions for a finite quandle to be Alexander. We then
describe an algorithm for determining finding all possible Alexander 
presentations of a finite quandle given its matrix and provide a URL for 
Maple code implementing this algorithm. 
\end{abstract}

\textsc{Keywords:} Alexander quandles, finite quandles, symbolic computation

\textsc{2000 MSC:} 57M27

\section{\large \textbf{Introduction}}

A \textit{quandle} is a set $Q$ with a binary operation $\tr:Q\times Q\to Q$
satisfying the three axioms
\newcounter{q}
\begin{list}{(\roman{q})}{\usecounter{q}}
\item{For all $a\in Q$, $a\tr a=a$,}
\item{For all $a,b \in Q$, there exists a unique $c\in Q$ such that 
$a=c\tr b$, and}
\item{For all $a,b,c\in Q$, we have $(a\tr b)\tr c=(a\tr c)\tr(b\tr c)$.}
\end{list}

The uniqueness of $c$ in axiom (ii) implies that the map $f_b:Q\to Q$
defined by $f_b(a)=a\tr b$ is a bijection for all $b\in Q$; we denote the
inverse $f^{-1}_b(a)$ by $a\tl b$. Then $Q$ forms a quandle under
the operation $\tl$, called the \textit{dual} of $(Q,\tr)$; in addition to
satisfying the analogs of the above axioms, $\tl$ also distributes over $\tr$ 
and vice-versa.

Quandles have been studied (indeed, rediscovered) numerous
times by various authors including Conway and Wraith, Brieskorn \cite{B},
Mateev \cite{Ma} and Fenn and Rourke \cite{FR}. The definition and notation
above were introduced by David Joyce in \cite{J}. 

Quandles and finite quandles in particular are of interest to knot theorists 
since associated to every knot there is a quandle, the \textit{knot quandle}, 
which is a complete invariant of knot type up to homeomorphism of topological 
pairs. Finite quandles then give us a convenient way to distinguish knots, 
since if two knot quandles $Q_1$ and $Q_2$ are isomorphic, the sets 
$\mathrm{Hom}(Q_1,Q)$ and $\mathrm{Hom}(Q_2,Q)$ must have the same number of 
elements. More sophisticated knot invariants involving counting homomorphisms 
to a finite quandle weighted by cocyles in various quandle cohomology 
theories are studied in various recent papers such as \cite{C1} and \cite{C2}.

One standard example of a quandle structure is the \textit{conjugation quandle}
of a group $G$. Specifically, $\mathrm{Conj}(G)$ has the same underlying
set as the group $G$ with quandle operation given by 
\[a\tr b=b^{-1}ab.\] 
Moreover, a subset of a group need not be a subgroup to form a quandle under 
conjugation; any union of conjugacy classes in a group forms a quandle. If 
the group or collection of conjugacy classes is finite, then we have a finite 
quandle. Other standard examples of finite quandles include the 
\textit{cyclic} quandle $Q=\{1,2,\dots, n\}$ with quandle operation defined
by
\[i\tr j = 2j-i \ (\mathrm{mod} \  n)\]
and the \textit{trivial 
quandle} $T_n=\{1,2,\dots, n\}$ with quandle operation given by
\[ i\tr j = i \quad \forall i,j\in T_n.\] The conjugation quandle of any 
abelian group is trivial.

Another example of a useful quandle structure is the \textit{Alexander quandle}
construction described in section \ref{AQ}. Alexander quandles have been 
studied in various papers (\cite{G}, \cite{N}, \cite{M}, \cite{MN1}), and most 
of the computations of quandle cohomology and counting invariants in recent 
papers have used finite Alexander quandles. In \cite{I}, the counting 
invariant is shown to depend on the classical Alexander invariants when the 
target quandle is Alexander, and in \cite{L} the quandle cohomology invariants 
for quadratic Alexander quandles are shown to be determined by the Alexander 
invariants for torus knots. Moreover, methods for computing the second
cohomology groups for Alexander quandles are given in \cite{M}, which permit
computation of the 2-cocyle invariants when the target quandle is Alexander.
Hence, when studying knots using quandle counting invariants, it is useful 
to know whether the target quandle is isomorphic to an Alexander quandle.

In \cite{HN}, a method was described for representing finite quandles as 
square matrices. These matrices can then be used to find all possible quandle 
structures of a given cardinality $n$. Specifically, for a quandle 
$Q=\{x_1,\dots, x_n\}$, the \textit{matrix of Q}, $M_Q$, is the matrix 
abstracted from the quandle operation table by dropping the $x$s and keeping 
only the subscripts. That is, $M_{ij}=k$ where $x_i\tr x_j=x_k$. This matrix 
notation was used in \cite{HN} to determine all quandle structures with up to 
5 elements. An improved algorithm for finding all quandle matrices 
together with URLs for the (rather large) files containing the results for 
$n=6,\ 7$ and 8 as well as a method for computing the counting
invariant using a target quandle given by a matrix are given in \cite{MN2}. 
An independently derived list of quandles of order less than or equal to six
can be found in \cite{C3}.

In this paper, we describe a method of determining whether a quandle
defined by a matrix is isomorphic to an Alexander quandle. In section 
\ref{AQ} we collect definitions, examples and necessary conditions for a 
finite quandle to be Alexander, as well as some results which are useful for 
the following section. We then describe an algorithm for taking 
a finite quandle matrix and finding all possible Alexander presentations of 
the given quandle, or determining when none exist. In section \ref{M} we 
describe our Maple implementation of this algorithm and provide a 
URL for this implementation.

\section{\large \textbf{Alexander quandles}} \label{AQ}

We begin with a definition.

\begin{definition} \textup{
Let $\Lambda=\mathbb{Z}[t^{\pm 1}]$ be the ring of Laurent polynomials
in one variable with integer coefficients. Let $M$ be a module over 
$\Lambda$. Then $M$ is a quandle, called an \textit{Alexander quandle}, 
with quandle operation given by \[ a\tr b = ta + (1-t)b.\]}
\end{definition}

\begin{example}\textup{
The trivial quandle $T_n$ is an Alexander quandle, namely the quotient 
module $T_n=\Lambda/(n,1-t)$:
\[ a\tr b = t(a)+(1-t)b =1(a) + (1-1)b = a.\]}
\end{example}

\begin{example}\textup{
Let $n\in \mathbb{Z}_+$ and $h\in \Lambda$. Then the quotient ring 
$\Lambda/(n,h)$ of Laurent polynomials modulo the ideal generated by
$n$ and $h$ is an Alexander quandle. More generally, an Alexander 
quandle may be a direct sum of such quotients or have a more 
complicated $\Lambda$-module structure. See \cite{N} for more examples.}
\end{example}

The structure of Alexander quandles has been explored in \cite{G} and 
\cite{N}. In particular, in \cite{N} we find

\begin{theorem}
If $M$ and $N$ are finite Alexander quandles, then there is an isomorphism
of Alexander quandles $\phi:M\to N$ iff there is an isomorphism of
$\Lambda$-modules $f:(1-t)M\to (1-t)N$.
\end{theorem}

See also \cite{AG} lemma 1.23.

This theorem tells us when two Alexander quandles are isomorphic, but how do we
know whether a quandle given, say, by a matrix, might be secretly Alexander? 

\begin{definition}\textup{
Let $Q$ be a finite quandle. The \textit{Alexanderization} of $Q$, denoted
$A(Q)$, is the free $\Lambda$-module on $Q$ modulo the submodule spanned
by elements of the form 
\[ tx_i+(1-t)x_j - x_i\tr x_j\]
for all $i,j=1,\dots,|Q|.$}
\end{definition}

Now suppose there is an isomorphism of quandles $\phi:Q\to M$ where $M$ is
a finite Alexander quandle. Then $\phi$ must factor through $A:Q\to A(Q)$, so
the diagram 
\[
\begin{array} {rcll}
Q & \longrightarrow           & A(Q)  \\
  & \searrow^{\phi} & \downarrow  \\
  &                           & M  
\end{array}
\]
must commute. In particular, injectivity of $\phi$ implies $A$ must also be 
injective. Thus we have

\begin{proposition} \label{nec} 
If a finite quandle $Q$ is isomorphic to an Alexander quandle, then 
$A:~Q\to~A(Q)$ is injective.
\end{proposition}

Proposition 2 gives us a way of identifying certain quandles as non-Alexander. 
For example, in the Alexanderization of the quandle with matrix
\[ Q = \left[ \begin{array}{lll}
 1 & 1 & 2 \\
 2 & 2 & 1 \\
 3 & 3 & 3
\end{array}\right]\]
we have $tx_1+(1-t)x_3 =x_1 \ \Rightarrow \ (1-t)x_1=(1-t)x_3$ from entry
$Q_{32}$. Then $Q_{13}$ says 
\[
x_2 =  tx_1+(1-t)x_3 =  tx_1 + (1-t)x_1 = x_1
\]
and $A$ is not injective; hence $Q$ is not Alexander.

\begin{definition} \textup{
A quandle $Q$ is \textit{abelian} if for all $a,b,c,d\in Q$ we have
\[(a\tr b) \tr (c\tr d) = (a\tr c)\tr(b\tr d).\]
}\end{definition}

\begin{proposition}
Alexander quandles are abelian.
\end{proposition}

\begin{proof}
Let $Q$ be Alexander. Then
\begin{eqnarray*}
(a\tr b) \tr (c\tr d) & = & t(ta+(1-t)b)+(1-t)(tc+(1-t)d) \\
 &=& t^2a +t(1-t)(b+c) +(1-t)^2d \\
& = & (a\tr c)\tr(b\tr d).
\end{eqnarray*}
\end{proof}

\begin{proposition}
If $Q$ is abelian, then $\tr$ is left-distributive. That is, 
\[a\tr(b\tr c)=(a\tr b) \tr (a\tr c).\]
\end{proposition}

\begin{proof}
Let $Q$ be an abelian quandle. Then for any $a,b,c\in Q,$ 
\[ (a\tr b) \tr (a\tr c) =  (a\tr a) \tr (b\tr c) =a\tr (b\tr c). \]
\end{proof}

\begin{corollary}
Alexander quandles are left-distributive. (See also \cite{LR}.)
\end{corollary}

These observations give us ways of testing whether a quandle is Alexander,
but they do not give any information about what the Alexander structure(s) 
on $Q$ might be. To solve this problem, we need a more constructive approach.

Let $M$ be an Alexander quandle. The facts that $t^{-1}\in \Lambda$ and 
$t(x+y)=tx+ty \ \forall x,y\in M$ show that multiplication by $t$ is an 
additive automorphism of $M$. Conversely, given any abelian group $A$ and 
automorphism $\phi\in \mathrm{Aut}_{\mathbb{Z}}(A)$, $A$ has the structure 
of a $\Lambda$-module, and hence an Alexander quandle, by defining 
$tx=\phi(x) \ \forall x\in M$. 

\begin{definition} \textup{
For any Alexander quandle $Q$, an \textit{Alexander presentation} 
consists of an abelian group structure on $Q$ together with an additive 
automorphism $\phi$ of this abelian group structure such that 
\[ a\tr b = \phi(a) +(1-\phi)(b)\quad \forall a,b\in Q, \]
that is, such that the induced Alexander quandle structure on $Q$ agrees 
with the original quandle structure.}
\end{definition}

\begin{definition}\textup{
Let $Q=\{x_1,x_2,\dots, x_n\}$ be a finite quandle. The \textit{(standard 
form) matrix of Q}, $M_Q$, is the matrix $M_Q$ whose entry in row
$i$ column $j$ is $k$ where $x_i\tr x_j=x_k$. A map $\phi:Q\to Q$ may be 
specified by a vector $v\in Q^n$ such that
\[v=\left[
\begin{array}{c}
\phi(1) \\
\phi(2) \\
\vdots \\
\phi(n)
\end{array}
\right].\]}
\end{definition}

\begin{example}\label{amex} \textup{
Let $Q$ be the Alexander quandle $\Lambda/(2,t^2+1)=\{x_1=0,x_2=1,x_3=t,
x_4=1+t\}$. Then $Q$ has quandle matrix
\[
\begin{array}{c|cccc}
    & x_1 & x_2 & x_3 & x_4 \\ \hline
x_1 & x_1 & x_4 & x_4 & x_1 \\
x_2 & x_3 & x_2 & x_2 & x_3 \\
x_3 & x_2 & x_3 & x_3 & x_2 \\
x_4 & x_4 & x_1 & x_1 & x_4 
\end{array}
\quad \longrightarrow
\quad M_Q=
\left[
\begin{array}{cccc}
1 & 4 & 4 & 1 \\
3 & 2 & 2 & 3 \\
2 & 3 & 3 & 2 \\
4 & 1 & 1 & 4 
\end{array}\right].
\]}
\end{example}

The matrix of a quandle is just the operation table of the quandle with
elements of the quandle replaced by their subscripts. Suppose we are given
a quandle matrix $M_Q$; we would like to determine whether $Q$ is an
Alexander quandle and, if it is, to find all Alexander presentations of $Q$.
Our basic method is to test the claim that $Q$ is Alexander by first checking
whether the Abelian condition is satisfied. If it is, we then proceed to 
try to construct all possible Alexander presentations of $Q$ from its quandle
structure.

We can represent finite abelian (or non-abelian) groups using a matrix 
notation very similar to our matrix notation for quandles.

\begin{definition} \textup{
Let $G=\{x_1,x_2,\dots, x_n\}$ be a finite group. The \textit{(standard 
form) Cayley matrix of G}, $C_G$, is the matrix $C_G$ whose entry in row
$i$ column $j$ is $k$ where $x_ix_j=x_k$ and $x_1$ is the identity 
element of $G$.}
\end{definition}

\begin{example} \textup{
The Alexander quandle in example 2 has abelian group structure
$\mathbb{Z}_2\times \mathbb{Z}_2$ with automorphism $\phi(a,b) = (b,a).$
An Alexander presentation for this quandle is}
\[C_Q=
\left[
\begin{array}{cccc}
1 & 2 & 3 & 4 \\
2 & 1 & 4 & 3 \\
3 & 4 & 1 & 2 \\
4 & 3 & 2 & 1 
\end{array}
\right], \quad \phi=
\left[
\begin{array}{c}
1 \\
3 \\
2 \\
4
\end{array}
\right].
\]

\textup{
We check that $\phi$ is an automorphism of $C_Q$ by checking that applying 
the permutation $\phi$ to each element of $C_Q$, then un-permuting the rows 
and columns by conjugating by the matrix of the permutation $\phi$ yields the
original matrix.}

\[
\left[
\begin{array}{cccc}
1 & 2 & 3 & 4 \\
2 & 1 & 4 & 3 \\
3 & 4 & 1 & 2 \\
4 & 3 & 2 & 1 
\end{array}
\right] =
\left[
\begin{array}{cccc}
1 & 0 & 0 & 0 \\
0 & 0 & 1 & 0 \\
0 & 1 & 0 & 0 \\
0 & 0 & 0 & 1 
\end{array}
\right]
\left[
\begin{array}{cccc}
1 & 3 & 2 & 4 \\
3 & 1 & 4 & 2 \\
2 & 4 & 1 & 3 \\
4 & 2 & 3 & 1 
\end{array}
\right]
\left[
\begin{array}{cccc}
1 & 0 & 0 & 0 \\
0 & 0 & 1 & 0 \\
0 & 1 & 0 & 0 \\
0 & 0 & 0 & 1 
\end{array}
\right].
\]
\end{example}

If $Q$ is a finite Alexander quandle, then one of the elements of $Q$
is the additive identity of $Q$ regarded as an abelian group. Then since
$a\tr 0 = ta +(1-t) 0 =ta$, the column corresponding to the additive identity
tells us the action of $t$ on $Q$. We can use this information to either 
recover the additive structure(s) of $Q$ or show that none is possible with 
the help of the following observations.

\begin{lemma} \label{ant}
If $M$ is an Alexander quandle, then for all $a,b\in M$ we have
\[a\tr b + b\tr a = a+b.\]
\end{lemma}

\begin{proof} If $M$ is Alexander, then for any $a,b\in M$
\begin{eqnarray*}
a\tr b + b\tr a & = & ta+(1-t)b + tb+ (1-t)a \\
 & = & ta + b -tb +tb +a -ta \\
 & = &  a+b.
\end{eqnarray*}
\end{proof}

\begin{lemma} \label{gab1}
Let $M$ be an Alexander quandle. Then for every $a,b,c\in M$ we have 
\[ a\tr b + b\tr c = a\tr c +b. \]
\end{lemma}

\begin{proof} 
If $M$ is Alexander, then for every $c\in M$, we have 
$a\tr c =ta +(1-t)c.$ Then
\begin{eqnarray*}
a\tr b + b\tr c & = & ta+(1-t)b +tb +(1-t)c \\
& = &  ta + b -tb +tb +(1-t)c \\
& = & a\tr c +b.
\end{eqnarray*}
\end{proof}

\begin{lemma} \label{gab2}
If $M$ is a finite Alexander quandle, then for a fixed $b\in M$ the map 
$g_b:M\to M$ defined by $g_b(x)=x+b$ is an automorphism of quandles 
(though not of modules). (See also \cite{AG}, \cite{LN}).
\end{lemma}

\begin{proof}
By definition,
\begin{eqnarray*}
g_b(x \tr y) & = &  x \tr y + b \\
& = & tx + (1-t)y + b \\
& = & tx + y - ty + b.
\end{eqnarray*}
On the other hand, 
\begin{eqnarray*}
g_b(x) \tr g_b(y) & = &  (x+b) \tr (y+b) \\
& = & t(x+b) + (1-t)(y+b) \\
& = & tx + tb + y + b -ty - tb \\
& = & tx + y + b - ty.
\end{eqnarray*}
and $g_b$ is a homomorphism of quandles. Setwise, $g_b$ is a cyclic 
permutation of the finite set $M=\{x_1,\dots, x_n\}$, and hence is
bijective. Thus, $g_b$ is a quandle automorphism of $M$.
\end{proof}

Applying lemmas \ref{gab1} and \ref{ant}, every entry of a finite quandle 
matrix tells us several equations that any Alexander structure on $Q$ must 
satisfy, each of which says that one entry in $C_Q$ is the same as another. 
In order to fill in the Cayley matrix, we need to have some starting values 
already filled in. If we assume that the additive identity element in $Q$ 
is $x_1$, then we can start the Cayley matrix with $C_Q[1,i]=C_Q[i,1]=i$ for 
each $i=1,\dots, n$. This assumption results in no loss of generality since 
by lemma \ref{gab2} if $Q$ is Alexander and $x_1$ is not the additive identity 
in $Q$, we simply apply the quandle automorphism $f(x_i)=x_i-x_1$ to obtain 
the same quandle operation table. With this assumption, multiplication by
$t$ in the Alexander structure on $Q$ is given by the first column of $Q$,
that is, $tx_i=x_{v[i]}$ where $v$ is the first column of $Q$.

Then for every element of the quandle matrix we compare the entries in $C_Q$
for each of the equations in lemmas \ref{ant} and \ref{gab1}; if the 
corresponding entries in $C_Q$ are different then $Q$ cannot be Alexander. If 
the entries are equal, or if both are blank, we move on; if one entry is known 
and the other blank, we fill in the blank with the known value. While going 
through this procedure, we exploit the facts that abelian groups are both 
commutative and associative to fill in the table and find contradictions more 
rapidly.

In this way, we fill in as much of the Cayley matrix as possible. There will 
generally be some entries which are left blank, since Alexander quandles may 
be isomorphic as quandles but distinct as $\Lambda$-modules. Thus, to find 
all possible Alexander structures on a given quandle, we systematically 
consider all possible ways of filling in the remaining blanks to obtain the
Cayley matrix of an abelian group. Having found all such matrices, it only 
remains to verify that the bijection given by the first column of the quandle 
matrix $Q$ is an automorphism of the abelian group structure so defined. 
If it is, then we have an Alexander presentation of the given quandle. 

Our implementation of this algorithm in Maple, available from 
\texttt{www.esotericka.org/quandles}, is described in the following 
section. 

\section{\large \textbf{Maple Implementation}} \label{M}

In this section we describe an implementation of the algorithm described in 
section \ref{AQ} in Maple. This code is available for download at 
\texttt{http://www.esotericka.org/quandles}; it uses the file 
\texttt{quandles-maple.txt} also available from the same website. Improvements
and bugfixes will be made as necessary.

We begin with some basic programs for working with abelian groups represented
by Cayley matrices. \texttt{assoctest} tests a matrix for associativity,
\texttt{commtest} tests a matrix for commutativity, and \texttt{invtest}
tests for the presence of inverses by checking that every row and column
contains the identity element 1. 

To implement the algorithm described in section \ref{AQ}, we start with 
\texttt{abgroupfill}, which uses the equations of lemmas \ref{ant} and 
\ref{gab1} to fill in entries in a standard form Cayley matrix, with zeroes 
representing unknown entries. The program compares the entries in the table 
which should be equal according to these equations, either replacing a zero 
with the nonzero value, doing nothing if both values are zero, or setting a 
``contradiction'' counter if it finds two different nonzero values, which 
results in quitting and reporting ``false''.

\texttt{abgroupfill} makes use of \texttt{cafill}, which uses associativity
and commutativity to fill in zeroes or find contradictions. The program runs 
through all triples, checking for associativity, then runs through all
pairs, checking for commutativity. When an entry is changed, a ``continue'' 
counter is set to true, so that the loop continues until no more zeroes
can be filled in, again exiting if a contradiction is found.

The program \texttt{zerofill} uses \texttt{findzero} to find the first zero 
entry in the matrix. It then fills in the zero with each possible nonzero 
entry from 1 to $n$, propagating these values through with \texttt{cafill}. 
Any resulting matrices are checked for rows or columns without the identity 
(which fail to be valid Cayley matrices); if every row and column contains 
either 1 or 0, then the matrix is added to the working list and the 
``continue'' counter is set. When no matrices in the working list have any 
zeroes, the loop is exited and we have a (possibly empty) list of Cayley 
matrices representing abelian groups. 

Our main program, \texttt{alextest}, first checks whether the input quandle is
abelian using \texttt{abqtest}. If it is, it then uses \texttt{abgroupfill} 
and \texttt{zerofill} to find all Cayley matrices corresponding to possible 
Alexander presentations of $Q$. The final step is to check whether the 
first column of $Q$ gives an automorphism of the Cayley matrix in question;
if it does, we have an Alexander presentation. We use the program 
\texttt{homtest} from the file \texttt{quandles-maple.txt} to check this.
If the list of Alexander presentations is empty, \texttt{alextest} returns 
``false.''

Next, we give some basic tools for constructing quandle matrices for
Alexander quandles. \texttt{cayley} returns the Cayley matrix for 
$\mathbb{Z}_n$. Many of the quandle tools in \texttt{quandles-maple.txt}
apply without modification to Cayley matrices; for example, \texttt{cprod}
gives the Cayley matrix of a cartesian product of two abelian groups,
\texttt{autlist} finds the automorphism group of an abelian group given
the Cayley matrix of the group, etc. \texttt{alexquandle} takes a Cayley 
matrix and a vector representing an automorphism of $C_Q$ and returns the
matrix of the resulting Alexander quandle structure. 

Finally, we include a short program \texttt{conjq} which takes a 
Cayley matrix for any group structure and returns the matrix of the 
conjugation quandle.

\end{document}